\documentclass[12pt]{article}
\usepackage{amsmath}

\title{Symmetry Classes of Alternating Sign Matrices} 
\author{David P. Robbins}
\date{August 4, 2000}

\def\TT{{\cal T}}      
\def\XX{{\cal X}}      
\def\ZZ{{\cal Z}}      
      
\def\d{\displaystyle}

\newtheorem{theorem}{Theorem}[section]
\newtheorem{conjecture}{Conjecture}[section]

\begin{document}

\maketitle

\begin{abstract}
An alternating sign matrix is a square matrix satisfying (i)   
all entries are equal to       
1, $-1$  or 0; (ii) every row and column has sum 1;    
(iii) in       
every row      
and column the non-zero entries alternate in sign.     
The 8-element group of symmetries of the square acts in an obvious way on      
square matrices.       
For any subgroup of the group of symmetries of the square      
we may consider the subset of matrices 
invariant under elements of this       
subgroup.  There are 8 conjugacy classes of these subgroups giving     
rise to 8 symmetry classes of matrices.
R. P. Stanley suggested the study of those alternating sign matrices in each of
these  
symmetry classes.  We have found  evidence suggesting that     
for six of the symmetry classes
there exist simple product formulas for the number of alternating sign 
matrices in the class. 
Moreover the factorizations of certain of their generating     
functions point to rather startling connections between several of the 
symmetry classes and cyclically symmetric plane partitions.    
       
\end{abstract}      
       
\newpage

\section{Introduction}

An alternating sign matrix is a square matrix satisfying       
\begin{itemize}       
\item[(i)] all entries are equal to 1, $-1$ or 0;     
\item[(ii)] all rows and columns have sum 1;  
\item[(iii)] in every row and column the non-zero entries alternate in sign.  
\end{itemize}
When we have occasion to consider the entries of an $n$ by $n$ alternating     
sign matrix, we shall always index these as $a_{i j }$ where $i$ and $j$ vary  
from 0 to $n-1$.       
       
Alternating sign matrices have been studied\footnote{I am posting this article at the request of Greg Kuperberg.
It it essentially the same manuscript that I wrote sometime in the late
1980's.  There has been much progress since that time particularly
by Kuperberg and Doron Zeilberger, but, except for indicating reference \cite{B}, I have made no attempt here
to bring the references up to date.}
in \cite{MRR1},
\cite{MRR2}, \cite{MRR3}, \cite{MRR4}, and \cite{RR}.  A much more recent and up to date
account is given in David Bressoud's book \cite{B}.  Alternating sign matrices have been
the source of a large number of conjectures related to their
enumeration.  In \cite{MRR2} and again in \cite{MRR3} there occurred
discussions of alternating sign matrices subject to certain symmetry
conditions.  R. P. Stanley suggested putting this study in a more
systematic framework in analogy with similar results described in
\cite{S} on symmetry classes of plane partitions.
       
The group of symmetries of the square acts in an obvious way on
square alternating sign matrices.      
For each subgroup of the group of symmetries of the square one 
may consider the set of alternating sign matrices invariant under this 
subgroup.  We call these sets the {\it symmetry classes}       
of alternating sign    
matrices.  This paper is primarily a list of conjectures related to these      
symmetry classes.      
       
There are 8 conjugacy classes of subgroups of the group of     
symmetries of the square.  These give  
rise to 8 symmetry classes of alternating sign matrices.       
Several of the symmetry classes have been previously studied,  
but we repeat the results here for completeness.       
We describe the symmetry classes next. 
\begin{center}
\begin{tabular}{cll}
\hbox{type} & condition &                                               \\   
1  &                               & \hbox{no conditions}               \\    
2  &  $ a_{i j }=a_{i , n - 1 - j } $             & \hbox{vertical axis}       
       \\     
3  &  $ a_{i j }=a_{n - 1 - i , n - 1 - j }$          & \hbox{half turn}       
           \\ 
4  &  $ a_{i j }=a_{j i } $                  & \hbox{diagonal} 
  \\  
5  &  $ a_{i j }=a_{j , n - 1 - i } $             & \hbox{quarter turn}
       \\     
6  &  $ a_{i j }=a_{i , n - 1 - j } =a_{n - 1 - i , j }$    & \hbox{horizontal 
and vertical}     \\  
7  &  $ a_{i j }=a_{j i }=a_{n - 1 - j , n - 1 - i }$      & \hbox{both
diagonals}              \\    
8  &  $ a_{i j }=a_{j i }=a_{i , n - 1 - j } $         & \hbox{all symmetries} 
\end{tabular}          
\end{center}
       
We have computed the number of alternating sign matrices in each of the
symmetry classes for small matrices.   
The numbers are shown in the first table of Section~\ref{sec:4}. 
       
In many cases  
we can conjecture      
simple formulas expressing the number of  $n$ by $n$ alternating       
sign matrices in a symmetry class as a product of small integers.      
These formulas are given in the second table of Section~\ref{sec:4}.

Our remaining results are concerned with       
generating functions for these 
symmetry classes.      
They   
point to rather startling connections between several of the   
symmetry classes and cyclically symmetric plane partitions.    
       
\section{Some Known Generating Functions}\label{sec:2}      
       
Many of
the generating functions for symmetry classes of alternating sign matrices     
can be expressed in terms of three     
functions, $Z_{n }(x,y,\mu )$, $T_{n }(x,\mu )$ and $R_{n }(x,\mu )$, which we 
describe below 
and which have been proved elsewhere to have some interesting properties.      
       
Let $\mu $ be a non-negative integer and define
$$Z_{n }(x,y,\mu )=\det(\delta _{i j }+z_{i j })_{0\le i,j<n}  $$      
where  
$$     
z_{i j }= \sum _{t,k=0}^{n-1}{\binom{i+\mu }{ t}}{\binom{k ^€}{ t}}{\binom{j-k+\mu -1 }{ 
j-k}}x^{k - t } 
$$     
for $0\le i<n-1$, $0\le j<n$ and       
$$ z_{n - 1 , j }= \sum _{t,k,l=0}^{n-1}{\binom{n-2+\mu -l}{ t-l}}
{\binom{k }{ t}}{\binom{j-k+\mu -1 }{ j-k}}x^{k - t }y^{l + 1 }  $$  
for $0\le j<n$.
       
We assign to $Z_{0 }(x,y,\mu )$ the conventional value 1.  When
regarded as a polynomial in $y$, its coefficients are known to form a
palindromic sequence.  The cases $\mu =0$ and $\mu =1$ will be most
important for us, so we list several of the values of $Z_{n }(x,y,0)$
and $Z_{n }(x,y,1)$ in the last table in Section~\ref{sec:4}.
       
Let    
$$     
T_{n }(x,\mu )=\det \left( \sum _{t = 0 }^{2 n - 2 }{\binom{i+\mu  }{ t-i}}       
{\binom{j }{ 2j-t}} x^{2 j - t }\right)_{0\le i,j<n}.  $$ 
The values of $T_{n }(x,0)$ and $T_{n }(x,1)$ will be particularly important for       
us.  Several are tabulated in the last section.
       
Let    
$$ Y(i,t,\mu )={\binom{i+\mu }{ 2i+1+\mu -t}}+{\binom{i+1+\mu  }{ 2i+1+\mu -t}} $$   
and define $R_{0 }(x,\mu )=1$ and      
$$ R_{n }(x,\mu )=     
\det \left(\sum _{t = 0 }^{2 n - 1 }Y(i,t,\mu )Y(j,t,0)x^{2 j + 1 - t  
}\right)_{i,j=0,\dots,n-1} $$ 
for $n\ge 1$.  
Some of the polynomials $R_{n }(x,0)$ and $R_{n }(x,1)$
are given in the last section of this paper.   
       
The first two  
functions are known to be the generating functions for classes   
of combinatorial objects.      
       
Let $\ZZ_{n }(\mu )$ be the set of shifted plane partitions (arrays of positive
integers)      
\[
\begin{matrix}
          a_{1 1 } & a_{1 2 }   & a_{1 3 } &\dots & \quad\quad  a_{1,\lambda  
_{1 }}\hfil   \\\      
              & a_{2 2 }   & a_{2 3 } &\dots & \quad a_{2,\lambda _{2 }}\hfil 
     \\     
                  &                 & \vdots & \cr   
              &       & a_{r r }       & õä& a_{r,\lambda _{r }}
\end{matrix}
\]
with $\lambda _{1 }>\cdots>\lambda _{r }$,     
with strictly decreasing columns and weakly decreasing 
rows, with no row length exceeding $n$ and     
such that the first entry $a_{i i }$   
of each row exceeds the length $\lambda _{i }-i+1$ of that     
row by precisely $2\mu $.  We say that a part $a_{i j }$ of such a plane       
partition      
is {\it special} if $\mu <a_{i j }\le j-i+\mu $.  To such a plane partition    
we assign a weight of $x^{r }y^{s }$ if there are $r$ special parts    
and $s$ parts in the first row equal to $n+2\mu $.     
Then one may verify,   
by arguments like those in \cite{MRR4}, that the sum of the   
weights of all partitions in $\ZZ_{n }(\mu )$  
is the polynomial      
$Z_{n }(x,y,\mu )$.    
       
We also remark that  $\ZZ_{n }(0)$ is known to be in one-to-one
correspondence with the class of cyclically symmetric plane partitions 
whose Ferrers graphs are contained in the box  
$$\XX_{n }=[1,n]\times [1,n]\times[1,n]. $$    
$\ZZ_{n }(1)$ is known to be in one-to-one correspondence      
with the class of all ``descending plane partitions'' with no parts    
exceeding $n+1$ (which are conjectured to be in one-to-one correspondence      
with $n+1$ by $n+1$ alternating sign matrices.)
       
Let $\TT_{n }(\mu )$ be the set of triangular arrays of positive integers      
\[
\begin{matrix}
          a_{1 1 } & a_{1 2 } & \dots   &\quad\quad a_{1 , n - 1 }   \\
          a_{2 1 } & a_{2 2 } & \dots   & a_{2 , n - 2 }   \\
          \vdots & \vdots && \\
          a_{n - 1 , 1 } &       &        & \\
\end{matrix}     
\]
such that      
all rows and columns are weakly decreasing     
and    
$a_{i 1 }\le n-i+1+\mu $ for $i=1,\dots,n-1$. 
We say that a part $a_{i j }$ of such a plane partition is {\it special} if    
$a_{i j }\le j$.       
If such a partition has $r$ special parts, we assign it a weight of $x^{r }$.  
Then one may verify that the generating function for $\TT_{n }(\mu )$ is       
precisely $T_{n }(x,\mu )$.    
       
We remark that $\TT_{n }(\mu )$ is known to be in one-to-one correspondence    
with the set of partitions in $\ZZ_{n }(\mu )$ which are invariant under a     
certain
involution.  In particular when $\mu =0$ these are in one-to-one correspondence
with cyclically symmetric plane partitions which are equal to their    
transpose-complements. 
(see \cite{MRR4} for definitions).  When $\mu =1$,    
these are in one-to-one correspondence with descending plane partitions
invariant under the involution described in \cite{MRR2}.   The argument is sketched   
in \cite{MRR4}.  The involution can be used to prove that in $Z_{n }(x,y,\mu )$       
has a palindromic coefficient sequence when regarded as a polynomial in $y$.

Some other properties of the polynomials       
$Z_{n }(x,y,\mu )$ and $T_{n }(x,\mu )$ are known.     
       
\begin{theorem}\label{thm:1}
Let $n$ be a non-negative integer.  Then  
\begin{align*} 
Z_{2 n }(x,1,\mu )&=T_{n }(x,\mu )R_{n }(x,\mu ) \\
\hbox{and}\quad  Z_{2 n + 1 }(x,1,\mu )&=2T_{n + 1 }(x,\mu )R_{n }(x,\mu ).    
\end{align*}
\end{theorem}

This is proved in \cite{MRR4}.
       
Below we use the abbreviation  
$$  (X)_{j }=X(X+1)(X+2)\cdots(X+j-1).  $$     
The hardest result is due to Andrews \cite{A}.  We repeat its statement here.       
       
\begin{theorem}[Andrews]\label{thm:2}
Let    
\begin{align*}    
\Delta _{0 }(\mu )&=2, \\     
\Delta _{2 j }(\mu )&={\frac{(\mu +2j+2)_{j }({\frac{\mu  }{ 2}}+2j+{\frac{3 }{ 2}})_{j - 1 }
}{  
(j)_{j }({\frac{\mu  }{ 2}}+j+{\frac{3}{ 2}})_{j - 1 }}},\qquad j>0, \\
\Delta _{2 j - 1 }(\mu )&={\frac{(\mu +2j)_{j - 1 }({\frac{\mu  }{ 2}}+2j+{\frac{1 }{ 2}})_{j
} }{
(j)_{j }({\frac{\mu  }{ 2}}+j+{\frac{1}{ 2}})_{j - 1 }}},\qquad\quad\
\  j>0. 
\end{align*}      
\end{theorem}

Then   
$$ Z_{n }(1,1,\mu )    
= \prod _{k = 0 }^{n - 1 }\Delta _{k }(2\mu ). $$      
       
The proof in \cite{A} is difficult.  An easier proof which applies only to  
$\mu =0$ and $\mu =1$ is given in \cite{MRR1}.

Using Theorems~\ref{thm:1} and \ref{thm:2}, one can prove  
a similar formula for $T_{n }(1,\mu )$.
       
\begin{theorem}\label{thm:3}   
$$     
T_{n }(1,\mu )=
{\frac{1 }{ 2^{n }}} \prod _{k = 0 }^{n - 1 }\Delta _{2 k }(2\mu ). $$     
\end{theorem}

A proof is given in \cite{MRR4}.      
       
Note that a consequence of Theorems~\ref{thm:2} and \ref{thm:3} is
that, for fixed integral $\mu $, when we substitute $x=1$ and $y=1$ in
$Z_{n }(x,y,\mu )$ or $T_{n }(x,\mu )$ we obtain a product of small
integers.
       
We also remark that implicit in \cite{MRR1} are formulas for  
$Z_{n }(1,y,0)$  and  $Z_{n }(1,y,1)$, but we have no need for these here.     
       
\section{Generating Functions for the Symmetry\\ Classes}\label{sec:3}
       
Now we proceed to describe what has been observed or proved concerning
the symmetry classes of alternating sign matrices.  There are eight
symmetry classes and we discuss these in Subsections~\ref{sec:3.1}
through \ref{sec:3.8}.
       
\subsection{All Alternating Sign Matrices}\label{sec:3.1}      
       
An alternating sign matrix always has a 1 in the top row.  Suppose that
$(a_{i j })$ is an alternating sign matrix and that $a_{0 s }=1$ and that the  
number 
of entries equal to $-1$ is precisely $r$.  Then we assign a weight of 
$x^{r }y^{s }$ to this matrix.  Let $A_{n }(x,y)$ be the ordinary generating   
function       
for all $n$ by $n$ alternating sign matrices, that is, 
the sum of the weights of all  
these matrices.
Then   
it is conjectured that 

\begin{conjecture}\label{conj:3.1.1.}    
$A_{n }(x,y)=Z_{n - 1 }(x,y,1)$.       
\end{conjecture}

\subsection{Flip Symmetric Alternating Sign Matrices}\label{sec:3.2} 
       
This symmetry class is empty unless    
$n=2m+1$ is odd.       
In this case there can be no zeros in column $m$ so that       
$a_{j m }=(-1)^{j }$.  If $k$ is the number of $-1$'s in the first $m$ columns 
of a flip symmetric alternating sign matrix, we assign it a weight     
of $x^{k }$.   Let $F_{n }(x)$ be the ordinary generating function  for this   
symmetry class of alternating sign matrices.

\begin{conjecture}\label{conj:3.2.1}    
$  F_{2 n + 1 }(x)=T_{n }(x,1).  $     
\end{conjecture}       

\subsection{Invariant Under the Half-Turn in its Own Plane}\label{sec:3.3}    

We assign to alternating sign matrices in this symmetry class a weight 
$x^{r }y^{s }$ where $a_{0 s }=1$ and  
$r$ is the number of orbits of $-1$'s under the action of      
the 2-element group generated by the   
half-turn.     
This number is half the number of $-1$'s in the
matrix unless $n=2m+1$ is odd and $a_{m m }=-1$ in which case the number of    
$-1$'s is an odd number
$2l+1$ and $r=l+1$.  We denote the generating function for  half-turn  
invariant $n$ by $n$ alternating sign matrices  by $H_{n }(x,y)$.      
       
For even size matrices we have observed the surprising 
formula
       
\begin{conjecture}\label{conj:3.3.1}    
$  H_{2 n }(x,y)=Z_{n }(x,y,0)Z_{n - 1 }(x,y,1). $     
       
$Z_{n }(x,y,0)$ is the generating function for cyclically  symmetric   
plane partitions with Ferrer's graph contained in the box $\XX_{n }$ and       
$Z_{n - 1 }(x,y,1)$ is conjectured to be the generating function for $n$ by $n$
alternating    
sign matrices. Thus    
Conjecture C1 suggests a one-to-one correspondence between     
half-turn symmetric $2n$ by $2n$ alternating sign matrices and the Cartesian   
product of the set of $n$ by $n$ alternating sign matrices and the     
cyclically symmetric plane partitions with Ferrers graph contained in  
the box $\XX_{n }$.    
\end{conjecture}
       
It is previously been observed that when we substitute $x=2$ in
generating functions for alternating sign matrices, they generally become      
much simpler.  For example, using the methods of       
\cite{MRR2} and \cite{RR} we can prove the following result.

\begin{theorem}\label{thm:4}   
\begin{align*}   
{\frac{H_{4 n }(2,1) }{ H_{4 n - 2 }(2,1)}}&=2^{2 n - 1 }{\frac{ {\binom{4n }{ 2n}}}{{\binom{2n 
}{ n}}  }} \\      
{\frac{H_{4 n + 2 }(2,1) }{ H_{4 n }(2,1)}}&=2^{2 n + 1 }{\frac{ {\binom{4n }{ 2n}}}{{\binom{2n 
}{ n}}  }} \\      
 H_{2 n + 1 }(2,1)&=2^{n }H_{2 n }(2,1).       
\end{align*}
\end{theorem}
       
For odd size half-turn symmetric alternating sign matrices the 
generating functions do not seem to factor.    
However when we set $y=1$, then $H_{2 n + 1 }(x,1)$ does seem to factor.       
In fact we have observed that  

\begin{conjecture}\label{conj:3.3.2}    
\begin{itemize}       
\item[(i)] $H_{4 n + 1 }(x,1)=R_{n }(x,0)T_{n }(x,1)S_{4 n + 1 }(x)$         
\item[(ii)] $H_{4 n - 1 }(x,1)=R_{n - 1 }(x,1)T_{n }(x,0)S_{4 n - 1 }(x)$     
where $S_{1 }(x)$,  $S_{3 }(x)$,$\dots$  are certain polynomials.     
\end{itemize}
       
The first few $S$'s are given in the table at the end of
Section~\ref{sec:4}.
\end{conjecture}
       
If we set $x=1$ in the generating function for odd size half-turn      
symmetric alternating sign matrices, obtaining the polynomials 
$H_{n }(1,y)$, 
we have no conjecture concerning       
the value of the resulting function.  Nevertheless these polynomials will      
appear later as factors of another generating function so we tabulate a few of 
their values in the last section.      
       
Its seems appropriate to repeat here   
an interesting observation from \cite{MRR3}.  
       
\begin{conjecture}\label{conj:3.3.3}
$  A_{n }(3,1)=3^{\hbox{deg}\;A_{n }(x,1)}H_{n }(1,1) .$       
\end{conjecture}       

\subsection{Equal to Transpose}\label{sec:3.4} 
       
We have enumerated the alternating sign matrices in this symmetry class
for $n=1,2,\dots,8$. Their numbers are
$1$, $2$, $5$,  $16$,  $67$, $2^{4 }\cdot 23$, $2\cdot 5\cdot263$,     
$2^{3 }\cdot 11\cdot 277$.  Apparently these numbers do not factor into small  
primes, so a simple product formula seems unlikely.  Of course this does       
not rule out other very simple formulas, but these would be more difficult     
to discover (let alone prove). 
       
\subsection{Invariant Under 90 Degree Rotation in its Own Plane}\label{sec:3.5}
       
When $n$ is even, each quarter of      
an alternating sign matrix in this symmetry class      
must have the sum of its entries       
equal to $n/4$.  It follows such matrices can exist for even $n$       
only when $n$ is actually a multiple of 4.     
Alternating sign matrices of all odd sizes do exist with this symmetry type.   
We assign to each such matrix a weight of      
$x^{r }y^{s }$ if $a_{0 s }=1$ and $r$ is the number of orbits of $-1$'s       
under the action of the group generated by the quarter turn.   
More precisely, $r$ is one-fourth the number of $-1$'s if the number of
$-1$'s is divisible by 4 or $l+1$ if the number of $-1$'s is $4l+1$    
(the only other possibility).  
It is clear that       
we must have $0<k<n-1$.  We denote by $Q_{n }(x,y)$ the ordinary generating    
function of $n$ by $n$ alternating sign matrices of this type. 
Then it follows that $Q_{n }(x,y)$ must be divisible by $y$.   
We have conjectures concerning $Q_{n }(x,1)$ and $Q_{n }(1,y)$.
       
\begin{conjecture}\label{conj:3.5.1} 
For $n\ge 1$       
\begin{align*}
 Q_{4 n }(1,y)&=yH_{2 n }(1,y,0)A_{n }(1,y)^{2 };\\
 Q_{4 n + 1 }(1,y)&=yH_{2 n + 1 }(1,y)A_{n }(1,y)^{2 }; \\
 Q_{4 n - 1 }(1,y)&=yH_{2 n - 1 }(1,y)A_{n }(1,y)^{2 }.    
\end{align*}    
\end{conjecture}
       
We have verified these formulas only for $n=1$, 2, 3 and 4.    
They suggest bijections between quarter-turn symmetric 
alternating sign matrices and various Cartesian products with  
other classes of alternating sign matrices and cyclically symmetric    
plane partitions.      
       
\begin{conjecture}\label{conj:3.5.2} 
There exists a sequence of polynomials $w_{0 }(x),w_{1 }(x),\dots$
such that $ Q_{2 n + 1 }(x,1)=w_{n }(x)w_{n + 1 }(x) $ if $n$ is even
and $Q_{2 n - 1 }(x,1)=xw_{n }(x)$ $w_{n + 1 }(x)$ if $n$ is odd.
Moreover $Q_{4 n }(x,1)=v_{n }(x)w_{2 n }(x)$ for suitable polynomials
$v_{n }(x)$.
\end{conjecture}       

The first statement has been observed to be true       
for  $n=0,\ldots,7$ and the second for $n=1,\dots,4$.      
Some of the $v$'s and $w$'s are tabulated in the last section.

We have a combinatorial interpretation for $w_{2 n }(x)$.  Consider
the cyclically symmetric plane partitions (described more fully in
\cite{MRR4}) with Ferrers graph $F$ contained in the box $\XX_{2 n
}=[1,2n]\times[1,2n]\times[1,2n] $ and which are self-complementary in
the sense that for $(i,j,k)\in \XX_{2 n }$, we have $(i,j,k)\in F$ if
and only if $(2n+1-i,2n+1-j,2n+1-k)\notin F$.  We say that a part
$a_{i j }$ of such a plane partition is {\it special} if $i\le a_{i j
}<j$.  Assign a weight to such a partition to be equal to $x^{k }$
where $k$ is the number of special parts.  Then the generating
function for this class of functions is equal to $x^{n }w_{2 n }(x)$
at least for $n\le 4$.  Also note that this is consistent with our
previous conjecture that the number of these self-complementary
cyclically symmetric plane partitions is the square of the number of
$n$ by $n$ alternating sign matrices.  (This was reported in \cite{S}.)
       
We note that if our other conjectures are true, then we would have     
$w_{1 }(1)$, $w_{3 }(1)$,$\dots$, equal the numbers of half-turn symmetric    
alternating sign matrices of sizes 1, 3, $\dots$.   Thus $w_{2 n + 1 }(x)$    
may be a generating function for these matrices with a suitable
assignment of weights.  However, we do not know of such a weight.  There       
is a similar observation connecting $v_{n }(x)$ with $2n$ by $2n$ half-turn    
symmetric alternating sign matrices.   
       
\subsection{Invariant Under Flips in Vertical and Horizontal
Axes}\label{sec:3.6}       
       
Alternating sign matrices with this symmetry can exist only when       
$n$ is odd.  When $n=2m+1$ we must have $a_{i m }=(-1)^{i }$ and $a_{m j       
}=(-1)^{j }$.  
       
For an alternating sign matrix of size $n=2m+1$ with this symmetry we  
assign a weight of $x^{k }$ where $k$ is the number of entries $a_{i j }=-1$   
with $0\le i,j<m$.  Let $P_{n }(x)$ be the generating function for the $n$ by  
$n$ alternating sign matrices satisfying this symmetry condition.  (We use     
$P$ since the axes of symmetry of the matrix look like a plus sign.)   
       
It is possible to conjecture a formula for the $P$'s terms of the $T$'s.       
       
\begin{conjecture}\label{conj:3.6.1} 
For $n\ge 1$       
\begin{align*}
P_{4 n + 1 }(x)&=T_{n }(x,1)T_{n }(x,0).  \\
P_{4 n - 1 }(x)&=T_{n - 1 }(x,1)T_{n }(x,0).
\end{align*}
This has been verified only for small $n$.     
\end{conjecture}       
       
\subsection{Invariant Under Flips in Both Diagonals}\label{sec:3.7}    
       
For $n=2,4,6,\dots,18$ we have found that the number of such
alternating sign matrices in this class is $2$, $2^{3 }$, $2^{2 }\cdot
13$, $8\cdot 71$, $2^{2 }\cdot 2609$, $2^{3 }\cdot 31 \cdot 1303$,
$2^{3 }\cdot 17 \cdot 124021 $, so apparently their number is not
given as a simple product of small integers.
       
On the other hand, when $n$ is odd, we can recognize a pattern, but
this is given in the second table of Section~\ref{sec:4}.  We have not
found any interesting properties of any of the generating functions.
       
\subsection{Invariant Under all Symmetries of the Square}\label{sec:3.8}
       
We have enumerated the alternating sign matrices in this symmetry class
and we find that for $n=1,3,5,\dots,17$       
the numbers of such $n$ by $n$ alternating sign matrices are   
1, 1, 1, 2, 4, 13, 46, $8\cdot 31$, $4\cdot 379$.      
Apparently these numbers do not factor into small      
primes, so again a simple product formula seems unlikely.

\section{Tables}\label{sec:4}       

\subsection{A Small Table of Numerical Values}\label{sec:4.1}   
       
In the table below the values not given have either not been computed  
or are too large to fit ($\ast$).   

\begin{tabular}{rrrrrrrrr}   
\multicolumn{9}{c}{Symmetry Type}\\    
\noalign{\medskip} 
\hbox{size}&1  &  2  &   3  &    4 &     5  &            6 &      7 &
8 \\   
\noalign{\medskip}  
   1&       1&       1&       1&       1&       1&       1&       1&       1 \\      
   2&       2&       0&       2&       2&       0&       0&       2&       0 \\       
   3&       7&       1&       3&       5&       1&       1&       3&       1 \\       
   4&      42&       0&      10&      16&       2&       0&       8&       0 \\       
   5&     429&       3&      25&      67&       3&       1&      15&       1 \\       
   6&    7436&       0&     140&     368&       0&       0&      52&       0 \\       
   7&  218348&      26&     588&    2630&      12&       2&     126&       2 \\       
   8&       $\ast$&       0&    5544&        &      40&       0&     568&       0 \\       
   9&       $\ast$&     646&   39204&        &     100&       6&    1782&       4 \\       
  10&       $\ast$&       0&  622908&        &       0&       0&   10436&       0 \\       
  11&       $\ast$&   45885& 7422987&        &    1225&      33&   42471&      13 \\       
  12&       $\ast$&       0&       $\ast$&        &    6460&       0&  323144&       0 \\       
  13&       $\ast$& 9304650&       $\ast$&        &   28812&     286& 1706562&      46 \\       
  14&       $\ast$&       0&       $\ast$&        &       0&       0&       $\ast$&       0 \\       
  15&       $\ast$&       $\ast$&       $\ast$&        & 1037232&    4420&       $\ast$&     248 \\       
  16&       $\ast$&       0&       $\ast$&        & 9779616&       0&       $\ast$&       0 \\       
  17&       $\ast$&       $\ast$&       $\ast$&       $\ast$&       $\ast$&  109820&       $\ast$&    1516        
\end{tabular}

\subsection{Table of Numerical Conjectures}\label{sec:4.2}      

Here the numbers of alternating sign matrices in classes       
1, 2, 3, 4, 6, and 7 are denoted by $A_{n }$, $F_{n }$, $H_{n }$, $Q_{n }$,    
$P_{n }$ and $X_{n }$. 
None of the results stated here have   
been proved.  The conjecture for the $P$'s is due to W.~H.~Mills. 
\begin{align*}
{\frac{A_{n + 1 } }{ A_{n }}} &= {\frac{\d{\binom{3n+1}{n}}}{\d{\binom{2n}{n}}}}.\\
\noalign{\medskip}     
{\frac{ F_{2n+1}}{F_{2n-1}}}&={\frac{\d{\binom{6n-2}{2n}}}{\d 2{\binom{4n-1      
}{ 2n}} }}. \\
\noalign{\medskip}     
{\frac{ H_{2 n + 1 }}{ H_{2 n }}}&= {\frac{\d{\binom{3n }{ n}}}{ \d{\binom{2n }{ n}} }};   
\qquad 
{\frac{ H_{2 n } }{ H_{2 n - 1 }}} = {\frac{\d 4{\binom{3n }{ n}}}{ \d 3{\binom{2n }{ n}} }}.\\  
\noalign{\medskip}     
Q_{4 n } &=H_{2 n }A_{n }^{2 };   \qquad       
   Q_{4 n + 1 }=H_{2 n + 1 }A_{n }^{2 }; \qquad
   Q_{4 n - 1 } =H_{2 n - 1 }A_{n }^{2 }. \\
\noalign{\medskip}     
  {\frac{ P_{4 n + 1 } }{ P_{4 n - 1 }}}   
&= {\frac{\d(3n-1){\binom{6n-3 }{ 2n-1}}}{ \d(4n-1){\binom{4n-2 }{ 2n-1}} }}; \qquad     
  {\frac{ P_{4 n + 3 } }{ P_{4 n + 1 }}}   
= {\frac{\d(3n+1){\binom{6n }{ 2n}}}{ \d(4n+1){\binom{4n }{ 2n}} }}. \\
\noalign{\medskip}     
{\frac{X_{2 n + 1 } }{ X_{2 n - 1 }}}&={\frac{ \displaystyle{\binom{3n }{ n}}       
}{ \displaystyle{\binom{2n-1 }{ n}} }}. 
\end{align*}      

\subsection{Tables of Generating Functions}\label{sec:4.3}     
\begin{align*}
Z_{1 }(x,y,0)&=1+y\\
Z_{2 }(x,y,0)&=2   + xy + 2y^{2 }\\
Z_{3 }(x,y,0)&=(4+x) + (4x+x^{2 })y + (4x+x^{2 })y^{2 }+ (4+x)y^{3 }\\
Z_{4 }(x,y,0)&=8+10x+2x^{2 }+(12x+15x^{2 }+3x^{3 })y   \\
         &\quad+(12x+15x^{2 }+4x^{3 }+x^{4 })y^{2 }+\dots\\
\noalign{\medskip}    
Z_{1 }(x,y,1)&=1+y\\
Z_{2 }(x,y,1)&=2 + (x+2)y + 2y^{2 }\\
Z_{3 }(x,y,1)&=(6+x) + (6+7x+x^{2 })y + (6+7x+x^{2 })y^{2 } + (6+x)y^{3 }\\
Z_{4 }(x,y,1)&=(24+16x+2x^{2 })+(24+52x+26x^{2 }+3x^{3 })y \\
&\quad+(24+64x+38x^{2 }+8x^{3 }+x^{4 })y^{2 }+\dots \\
\noalign{\medskip}     
      T_{1 }(x,0)&=1 \\
      T_{2 }(x,0)&=1+x \\
      T_{3 }(x,0)&=1+5x+4x^{2 }+x^{3 } \\
      T_{4 }(x,0)&=1+14x+49x^{2 }+62x^{3 }+34x^{4 }+9x^{5 }+x^{6 }\\    
\noalign{\medskip}     
      T_{1 }(x,1)&=1\\
      T_{2 }(x,1)&=2+x\\
      T_{3 }(x,1)&=6+13x+6x^{2 }+x^{3 }\\
      T_{4 }(x,1)&=24+136x+234x^{2 }+176x^{3 }+63x^{4 }+12x^{5 }+x^{6 }\\
\end{align*}      

\begin{align*}
      R_{1 }(x,0)&=4+x  \\
      R_{2 }(x,0)&=16+40x+9x^{2 }+x^{3 } \\
      R_{3 }(x,0)&=64+560x+1036x^{2 }+629x^{3 }+125x^{4 }+16x^{5 }+x^{6 } \\
\noalign{\medskip}     
      R_{1 }(x,1)&=6+x \\
      R_{2 }(x,1)&=60+70x+12x^{2 }+x^{3 } \\
      R_{3 }(x,1)&=840+3080x+3038x^{2 }+1224x^{3 }+195x^{4 }+20x^{5 }+x^{6 }\\  
\noalign{\medskip}     
H_{1 }(1,y)&=1 \\
H_{3 }(1,y)&=1+y+y^{2 } \\ 
H_{5 }(1,y)&=3+6y+7y^{2 }+6y^{3 }+3y^{4 } \\
H_{7 }(1,y)&=25+75y+123y^{2 }+142y^{3 }+123y^{4 }+75y^{5 }+25y^{6 }\\   
\noalign{\medskip}     
S_{1 }(x)&=1 \\
S_{3 }(x)&=2+x \\
S_{5 }(x)&=2+3x \\
S_{7 }(x)&=8+26x+7x^{2 }+x^{3 } \\
S_{9 }(x)&=12+74x+78x^{2 }+31x^{3 }+3x^{4 }    
\end{align*}      

\begin{align*}       
      w_{0 }(x)&=1                                                  \\
      w_{1 }(x)&=1                                                  \\
      w_{2 }(x)&=1                                                  \\
      w_{3 }(x)&=2+x                                                \\
      w_{4 }(x)&=3+x                                                \\
      w_{5 }(x)&=4+14x+6x^{2 }+x^{3 }                                       \\
      w_{6 }(x)&=15+25x+8x^{2 }+x^{3 }                                      \\
      w_{7 }(x)&=8+88x+222x^{2 }+192x^{3 }+65x^{4 }+12x^{5 }+x^{6 } \\
      w_{8 }(x)&=105+490x+665x^{2 }+386x^{3 }+102x^{4 }+15x^{5 }+x^{6 }\\   
\noalign{\medskip}     
      v_{1 }(x)&=2                                                  \\
      v_{2 }(x)&=4+6x                                               \\
      v_{3 }(x)&=8+52x+60x^{2 }+20x^{3 }                                    \\
      v_{4 }(x)&=16+272x+1212x^{2 }+2000x^{3 }+1470x^{4 }+504x^{5
      }+70x^{6 }
\end{align*}

\end{document}